\documentclass[12pt, reqno]{amsart}

\usepackage{amssymb,upref,enumerate}
\usepackage{tikz}
\usepackage{mathabx}
\usepackage{amsmath,amssymb,yfonts,amscd,amsthm,amsxtra,mathrsfs}
\usepackage{latexsym}
\usepackage{soul}
\usepackage{verbatim}
\usepackage[T1]{fontenc}
\usepackage{tgtermes}
\usepackage[latin9]{inputenc}

\usepackage{cite} 

\usepackage[
colorlinks=true,
linkcolor=blue,
citecolor=blue,
urlcolor=blue,
]{hyperref}

\makeatletter
\def\Ddots{\mathinner{\mkern1mu\raise\p@
\vbox{\kern7\p@\hbox{.}}\mkern2mu
\raise4\p@\hbox{.}\mkern2mu\raise7\p@\hbox{.}\mkern1mu}}
\makeatother

\def\Xint#1{\mathchoice
{\XXint\displaystyle\textstyle{#1}}%
{\XXint\textstyle\scriptstyle{#1}}%
{\XXint\scriptstyle\scriptscriptstyle{#1}}%
{\XXint\scriptscriptstyle\scriptscriptstyle{#1}}%
\!\int}
\def\XXint#1#2#3{{\setbox0=\hbox{$#1{#2#3}{\int}$}
\vcenter{\hbox{$#2#3$}}\kern-.5\wd0}}

\def\dashint{\Xint-}

\newtheorem{theorem}{Theorem}[section]
\newtheorem{nonum}{Theorem}

\newtheorem{conjecture}[theorem]{Conjecture}
\newtheorem{corollary}[theorem]{Corollary}

\theoremstyle{definition}

\newcommand{\absval}[1]{\mbox{$|#1|$}}
\newcommand{\norm}[1]{\mbox{$\left\| #1 \right\|$}}

\def\ep{{\epsilon}}

\def\al{{\alpha}}
\def\R{\mathbb R}

\def\ra{\rightarrow}

\def\bey{\begin{eqnarray*}}
\def\eey{\end{eqnarray*}}

\allowdisplaybreaks

\begin{document}


\

\title[]{A homage to Guido Weiss and his leadership of the Saint Louis team: Commutators of Singular Integrals and Sobolev inequalities}

\author{Cong Hoang}

\address{Cong Hoang \\
 Department of Mathematics \\
Florida Agricultural and Mechanical University \\
 Tallahassee, FL 32307, USA}

\email{cong.hoang@famu.edu}

\author{Kabe Moen}

\address{Kabe Moen \\
 Department of Mathematics \\
 University of Alabama \\
 Tuscaloosa, AL 35487, USA}

\email{kabe.moen@ua.edu}

\author{Carlos P\'erez}

\address{Carlos P\'erez \\
 Department of Mathematics \\
 University of the Basque Country \\
 Ikerbasque and BCAM \\
 Bilbao, Spain}

\email{cperez@bcamath.org}


\thanks{
K.M. is supported by Simons Collaboration Grant for Mathematicians, 160427.
C. P. is supported by grant  PID2020-113156GB-I00, Spanish Government; by the Basque Government through grant { IT1615-22 and the BERC programme 2022-2025 program  and by BCAM Severo Ochoa accreditation CEX2021-001142-S, Spanish Government.} He is also very grateful to the Mittag-Leffler Institute where part of this research was carried out. }

\maketitle

\vspace{-1cm}

\begin{abstract} We extend some classical Sobolev type inequalities for linear and non-linear commutators.  
\end{abstract}


\section{Introduction}

It is with great pleasure that we describe some of the results concerning commutators in \cite{CRW}, which has been one of the main mathematical contributions of Guido L. Weiss as part of the celebrated "Saint Louis team" during the 7th decade  of the XX century, in connection with some new Sobolev type estimates in the spirit of our recent work \cite{HMP}. This will be {outlined} in Section \ref{linear commutators}.
Similarly, in Section \ref{nonlinear commutators} we will {outline Sobolev type estimates as in \cite{HMP} for certain nonlinear commutators introduced in another relevant paper \cite{RW}. Another very influential work of Guido is  \cite{StW} which can be seen as model case of the work of Muckenhoupt and Wheeden \cite{MW} central in the applications of our results in \cite{HMP}.}

{We will begin with a perspective from our recent} work  \cite{HMP} where the classical Gagliardo-Nirenberg-Sobolev (GNS) inequalities have been generalized and extended in several directions.

\section{Recent extensions of the classical GNS inequalities} \label{Extension Sobolev}

\subsection{Historical remarks}

The classical Sobolev {inequalities from the 30's of the last century which became} the  GNS inequalities in the 60's of the same century have been improved in many directions. {The purpose of this paper is to survey some recent improvements and extensions obtained by the authors} in \cite{HMP}.  We begin by outlining the history of the GNS estimates in its different incarnations related to the context of $L^p$ spaces and in $\R^n$ for $n\geq 2$ 
{(see the recent treatment \cite{KLV} Section 3.2.). }

 \subsubsection{{\color{red}The origins}}

\

Our starting point is commonly referred to as the GNS inequalities which are ubiquitous in harmonic analysis and P.D.E. Hereafter, we will use the notation of $p^*=\frac{np}{n-p}$ for the Sobolev exponent, often expressed in the following way
$$\frac1p-\frac1{p^*}=\frac1n.$$

\begin{nonum} (The classical GNS). { Let $1\leq p<n$.  Then }
\begin{equation}\label{GNSineq} \|f\|_{L^{p^*}(\R^n)}\leq c\, \|\nabla f\|_{L^p(\R^n)}, \quad f\in C^\infty_c(\R^n).
\end{equation}

\end{nonum}

One common approach to proving \eqref{GNSineq}, typically only when $1<p<n$,  is via the use of the {well-known} pointwise estimate (see \cite{Sal, St, KLV})
\begin{equation}\label{represent} 
|f(x)|\leq c \,I_1(|\nabla f|)(x)
\end{equation} 
combined with the boundedness (a {well-known} classical result)
\begin{equation}\label{bounded} 
I_1: L^p(\R^n) \to L^{p^*}(\R^n)
\end{equation} 
%
for $p\in(1,n)$. $I_1$ is the Riesz potential operator of order $\al=1$ given by 
$$I_\al f(x)=\int_{\R^n}\frac{f(y)}{|x-y|^{n-\al}}\,dy, \quad 0<\al<n.$$

At the endpoint $p=1$ and {$1^*=\frac{n}{n-1}=n'$,} $I_1$ only satisfies the weak endpoint boundedness
$$I_1:L^1(\R^n)\ra L^{n',\infty}(\R^n),$$
yet the strong inequality \eqref{GNSineq} still holds in this case (more on this below).

\

\subsubsection{\color{red} The 70's, the influence of Muckenhoupt and Wheeden}

\

A relevant extension of \eqref{bounded} and hence of \eqref{GNSineq} in the case $p>1$ was obtained by B. Muckenhoupt and R. Wheeden \cite{MW}. 
They showed in this paper that $I_1$ satisfies weighted bounds of the form 
\begin{equation}\label{fracbdd}\|wI_1 f\|_{L^{p^*}(\R^n)}\leq c\,\|wf\|_{L^p(\R^n)} \end{equation}
 if and only if $w\in A_{p,p^*}$ {  namely}
$$[w]_{A_{p,p^*}}=\sup_Q \left(\dashint_Q w^{p^*} \right) \left(\dashint_Q w^{-p'} \right)^{\frac{p^{*}}{p'}}<\infty.$$
Again, at the endpoint $p=1$ only the weak boundedness holds, and Muckenhoupt and Wheeden proved that
\begin{equation}\label{weakbdd} \|I_1f\|_{L^{n',\infty}(w^{n'})}\leq c\,\|wf\|_{L^1(\R^n)} \end{equation}
holds if and only if $w\in A_{1,n'}$, i.e.
$$\left(\dashint_Q w^{n'}\right)\leq c\inf_Q (w^{n'})$$
where the smallest constant $c$ is denoted as $[w]_{A_{1,n'}}$. 

Combining  \eqref{represent} and \eqref{fracbdd} yields the following weighted GNS inequality.

\begin{nonum} (The weighted GNS) {  Let $1\leq p<n$ and let $w\in A_{p,p^*}$. Then}
\begin{equation}\label{GNSweight}
\|wf\|_{L^{p^*}(\R^n)}\leq c\,\|w\nabla f\|_{L^p(\R^n)}, \qquad f\in C^\infty_c(\R^n).
\end{equation}
\end{nonum}
\noindent When $p=1$ and $w\in A_{1,n'}$, the operator $I_1$ only satisfies a weak-type inequality, yet the strong bound holds for the weighted GNS as we will detail in the next section.  
%
%
\

\subsubsection{\color{red} The influence of modern weighted theory in harmonic analysis}
 
\

Besides its dependence on the dimension and $p$, the constant $c$ in \eqref{GNSweight} also depends on the constant $[w]_{A_{p,p^*}}$, and this relationship is of interest. 

The sharp weighted bound for $I_1$ was proven in \cite{LMPT}  
\begin{equation}\label{StrongI1LMPT}
\|wI_1 f\|_{L^{p^*}(\R^n)}\leq c_{p}\,[w]_{A_{p,p^*}}^{\frac1{n'}\max\{1,\frac{p'}{p^*}\}}\,\|wf\|_{L^p(\R^n)}, 1<p<n,
\end{equation}
which leads to the following quantitative weighted GNS inequality,
\begin{equation}\label{sobolevbadbound}
\|wf\|_{L^{p^*}(\R^n)}\leq c_{p}\,[w]_{A_{p,p^*}}^{\frac1{n'}\max\{1,\frac{p'}{p^*}\}}\,\|w\nabla f\|_{L^p(\R^n)}, \quad 1<p<n.
\end{equation}
{Inequality \eqref{sobolevbadbound} can be improved and  the sharp version of the GNS inequalities is contained in the next theorem which can also be found in \cite{LMPT}.
\begin{nonum} (The precise weighted GNS) {  Let $1\leq p<n$ and let $w\in A_{p,p^*}$. Then}
\begin{equation}\label{oneweightSob}
\|wf\|_{L^{p^*}(\R^n)}\leq c_{p}\,[w]_{A_{p,p^*}}^{\frac1{n'}}\,\|w\nabla f\|_{L^p(\R^n)}, \quad 1\leq p<n.
\end{equation}
\end{nonum}

Notice that inequality \eqref{oneweightSob} improves \eqref{sobolevbadbound} 
{in the exponent} and holds at the endpoint $p=1$}. In fact, it is the endpoint $p=1$ that is crucial in proving \eqref{oneweightSob}.
We now briefly explain how to obtain the sharp inequality \eqref{oneweightSob}
by using the theory of extrapolation of weights. {First, using the pointwise inequality \eqref{represent}, we claim the weak type endpoint inequality for $p=1$ and $w\in A_{1,n'}$, }
%
\begin{equation*}\label{endpoint} 
\|f\|_{L^{n',\infty}(w^{n'})}\leq c\,[w]_{A_{1,n'}}^{\frac1{n'}}\,\|w\nabla f\|_{L^1(\R^n)}.
\end{equation*}

Given any measure $\mu$, the quantity $\|g\|_{L^{n',\infty}(\mu)}$ is equivalent to a norm since $n'>1$. 
Since $I_1$ is a convolution operator with the kernel $k(x)=|x|^{1-n}$, we may use the Marcinkiewicz integral inequality for normed spaces to obtain
\begin{align*}
\|I_1 f\|_{L^{n',\infty}(\mu)}=\|k*f\|_{L^{n',\infty}(\mu)} & \leq c\int_{\R^n}|f(y)|\big\||\cdot-y|^{1-n}\big\|_{L^{n',\infty}(\mu)}\,dy \\
& \leq c\int_{\R^n}|f(y)| M\mu(y)^{\frac{1}{n'}}\,dy \end{align*}
{where we have used the the fact that $\big\||\cdot-y|^{1-n}\big\|_{L^{n',\infty}(\mu)}\approx M\mu(y)^{\frac1{n'}},$ with $M$ being the usual Hardy-Littlewood maximal function.
%
Thus $I_1$ satisfies the weak bound
\begin{equation}\label{fracbddgoodboundsweak}
\|I_1 f\|_{L^{n',\infty}(w^{n'})}\leq c\int_{\R^n}|f|M(w^{n'})^{\frac1{n'}}\leq  c\,[w]_{A_{1,n'}}^{\frac1{n'}}\,\|wf\|_{L^1(\R^n)}.
\end{equation} 
Inequality \eqref{fracbddgoodboundsweak} now implies 
\begin{equation}\label{fracLpbddgoodweak}
\|I_1 f\|_{L^{p^*,\infty}(w^{p^*})}\leq  c_{p}\,[w]_{A_{p,p^*}}^{\frac1{n'}}\,\|wf\|_{L^p(\R^n)}, \quad 1\leq p<n,
\end{equation} 
by the sharp extrapolation theorem from \cite[Theorem 2.1]{LMPT} which is based on the off-diagonal extrapolation theorem by
Harboure, Mac\'ias, and Segovia \cite{HMS}. Inequality \eqref{fracLpbddgoodweak} now implies the weak Sobolev inequality
\begin{equation}\label{weakSobolev} \|f\|_{L^{p^*,\infty}(w^{p^*})}\leq c\|I_1(|\nabla f|)\|_{L^{p^*,\infty}(w^{p^*})}\leq c_{p}[w]_{A_{p,p^*}}^{\frac1{n'}}\|w\nabla f\|_{L^p(\R^n)}.\end{equation} 
The truncation method of Maz'ya \cite{Maz} (see also the survey article \cite{Haj}) applied to \eqref{weakSobolev} now implies the strong inequality with the same dependence on $[w]_{A_{p,p^*}}$, namely, inequality \eqref{oneweightSob}.

The sharpness of \eqref{oneweightSob} follows from the following example. Let
$$f_\delta(x)=\exp(-|x|^\delta) \quad \text{and} \quad w_\delta(x)=|x|^{\frac{\delta-n}{p^*}}$$ 
where $\delta \in (0,1)$. Then, 
$$\|w_\delta f_\delta\|_{L^{p^*}(\R^n)} \approx \delta^{-\frac{1}{p^*}}, [w_\delta]_{A_{p,p^*}}\approx \delta^{-1}, \ \ \text{and} \ \|w_\delta\nabla f_\delta\|_{L^p(\R^n)}\approx \delta^{\frac1{p'}}$$
showing that {  the power} $\frac1{n'}$ in \eqref{oneweightSob} cannot be improved by taking $\delta\ra 0^+$.

%
%

\subsection{{The toy motivation:} the maximal function case} \label{ToyExample}
We are interested in extending inequality \eqref{oneweightSob} to an inequality with an operator on the left hand side. We begin with the Hardy-Littlewood maximal function defined by 
$$Mf(x)=\sup_{Q\ni x} \dashint_Q |f(y)|\,dy.$$
The classical $A_p$ theorem of Muckenhoupt showed that
$$M:L^p(w)\ra L^p(w)$$
if and only if  $w\in A_p${, i.e.}
$$[w]_{A_p}=\sup_Q\left(\dashint_Q w\right)\left(\dashint_Q w^{1-p'}\right)^{p-1}<\infty.$$

There is an intimate connection between the $A_p$ and $A_{p,p^*}$ classes which will be very useful in what follows, {namely}
$$
w \in A_{p,p^*} \ \ \mbox{if and only if} \ \ w^{p^*} \in A_{\frac{p^*}{n'}}
$$
and
$$ [w]_{A_{p,p^*}}=[w^{p^*}]_{A_{\frac{p^*}{n'}} }.$$
Hence, if $w\in A_{p,p^*}$, then the weight $w^{p^*}$ belongs to $A_{p^*}$ since $A_{\frac{p^*}{n'}}\subset A_{p^*}$.

{By Muckenhoupt's theorem we obtain the following Sobolev inequality
\begin{equation}\label{maxweightSob}\|wMf\|_{L^{p^*}(\R^n)}\leq c\, \|w\nabla f\|_{L^p(\R^n)},\end{equation}
 which is an extension of \eqref{GNSweight} since $Mf\geq |f|$, a.e. by the Lebesgue differentiation theorem.  }

Consider the problem of finding the dependence on $[w]_{A_{p,p^*}}$ of the constant $c$ in inequality \eqref{maxweightSob}.  {One approach is to use Buckley's sharp bound  in \cite{Buc} which yields 
\begin{equation}\label{BuckBd}\|Mf\|_{L^p(w)}\leq c p'\,[w]_{A_p}^{\frac{1}{p-1}}\|f\|_{L^p(w)}.\end{equation}
This inequality leads to the following calculation, 
\begin{multline*}\|wMf\|_{L^{p^*}(\R^n)} = \|Mf\|_{L^{p^*}(w^{p^*})}\leq c_{p}\,[w^{p^{*}}]_{A_{p^*}}^{\frac1{p^*-1}}\,\|f\|_{L^{p^*}(w^{p^*})} \\
\leq  c_{p}\,[w^{p^{*}}]_{A_{p^*}}^{\frac1{p^*-1}}\,   [w]_{A_{p,p^*}}^{\frac1{n'}}\,\|w\nabla f\|_{L^p(\R^n)}{\leq \,} c_{p}\,[w]_{A_{p,p^*}}^{\frac{1}{p^*-1}+\frac1{n'}}\|w\nabla f\|_{L^p(\R^n)},
\end{multline*}
using \eqref{oneweightSob}. We may improve these estimates by using the fact that the weight $w^{p^*}$ belongs to the better class $A_{\frac{p^*}{n'}}$. Here we will use the following observation by J. Duoandikoetxea \cite{Duo} (see also \cite{LPR}) who showed that if $w\in A_q$ for $1\leq q<p$, namely a stronger condition, then 
\begin{equation*}\label{subclassest}\|Mf\|_{L^p(w)}\leq  c_{p,q}\,[w]_{A_q}^{\frac1p}\|f\|_{L^p(w)}.\end{equation*}
This leads to 
\begin{equation}\label{betterbd}\|wMf\|_{L^{p^*}(\R^n)}\leq c\,[w]_{A_{p,p^*}}^{\frac1{p^*}+\frac1{n'}}\|w\nabla f\|_{L^p(\R^n)},
\end{equation}
using again \eqref{oneweightSob}. Observe that \eqref{betterbd} also holds for $p=1$.  
We do not know if the exponent is sharp as in the case of  \eqref{oneweightSob}.

A similar argument works for iterations of $M$,  namely
\begin{equation}\label{maxiteratedSob}
\|wM^m f\|_{L^{p^*}(\R^n)}\leq c_{m,p}\,[w]_{A_{p,p^*}}^{\frac{m}{p^*}+\frac1{n'}}\,\|w\nabla f\|_{L^p(\R^n)}, \quad 1\leq p<n,
\end{equation} 
where $M^{m}$ denotes the $m$-th  iteration of the
Hardy-Littlewood maximal operator, namely $M^{m}= \underbrace{M\circ \cdots  \circ  M}_{(m\, times)}.$

Recently, the estimates above have been significantly improved in \cite{HMP} by using a different approach.  In this work, it is shown that the bound}
\begin{equation} \label{maxSob}
\|wMf\|_{L^{p^*}(\R^n)}\leq c_{p}\,[w]_{A_{p,p^*}}^{\frac1{n'}\max\{1,\frac{p'}{p^*}\}}\,\|w\nabla f\|_{L^p(\R^n)}, 
\end{equation} 
holds when $p\in(1,n)$ although the result is worse near $p=1$  than \eqref{betterbd}.
%
This result follows applying $M$ to  \eqref{represent} using the fact that $I_1f \in A_1$ with dimensional constant, independent of $f$, namely
$$Mf(x)\leq c\,I_1(|\nabla f|)(x), \quad f\in C_c^\infty(\R^n).$$

The same pointwise bound holds for the iterated maximal operator, namely
$$M^mf(x)\leq c_{m}\,I_1(|\nabla f|)(x), \quad f\in C_c^\infty(\R^n),$$
and hence
\begin{equation}
\label{maxSobIter}
\|wM^m f \|_{L^{p^*}(\R^n)}\leq c_{m,p}\,[w]_{A_{p,p^*}}^{\frac1{n'}\max\{1,\frac{p'}{p^*}\}}\,\|w\nabla f\|_{L^p(\R^n)}, \quad  1<p<n.
\end{equation} 
Observe that \eqref{maxSobIter} is an extension of \eqref{oneweightSob} since, again, $M^mf\geq |f|$, a.e.  Furthermore, the exponent is better than \eqref{maxiteratedSob} since it does not depend on $m$.

When $p=1$, we cannot use the truncation method mentioned earlier to derive the sharp estimate \eqref{oneweightSob} since 
$M$ and $M^m$ are non-local operators.  Nevertheless, \eqref{maxSobIter} is also true, but with a worse constant, {depending on $m$,}
namely 
\begin{equation}
\|wM^mf\|_{L^{n'}(\R^n)} \leq c_{m}\,   [w]_{A_{1,n'}}^{\frac{m+1}{n'}}\,\|w\nabla f\|_{L^1(\R^n)}.
\end{equation} 
%

%
%
%
%

%

%
%
%
%

%

It is possible to further improve \eqref{maxiteratedSob} and hence \eqref{oneweightSob}. Indeed, if we define 
$$
M_r{f}= \big[M(|f|^r)\big]^{\frac1r} \quad \text{for }\, r\in(0,\infty),
$$
then it is possible to prove in the range $p\in(1,n)$ that
\begin{equation} \label{Sobimpro}
\|wM_{n'}f\|_{L^{p^*}(\R^n)}\leq c_{p}\,[w]_{A_{p,p^*}}^{\frac1{n'}\max\{1,\frac{p'}{p^*}\}}\,\|w\nabla f\|_{{L^p(\R^n)}},
\end{equation}
by showing the pointwise estimate (see \cite{HMP})
\begin{equation}\label{PointMn'}
M_{n'}f(x) \leq c\, I_{1}(|\nabla f|)(x),
\end{equation}
and using \eqref{StrongI1LMPT}.  Again, \eqref{Sobimpro} improves the estimate that is obtained by using Buckley's inequality \eqref{BuckBd}, namely the less desirable bound,
\begin{equation*} 
\|wM_{n'}f\|_{L^{p^*}(\R^n)}\leq c_p\,[w]_{A_{p,p^*}}^{  \frac1{p^*-n'} +\frac1{n'}}\|w\nabla f\|_{{L^p(\R^n)}}.
\end{equation*}
Moreover, we cannot use {$A_{\frac{p^*}{n'}}\subset A_{p^*}$}  because $A_{\frac{p^*}{n'}}$ is the precise class of weights for which $M_{n'}$ is bounded. 

{Observe that inequality \eqref{Sobimpro} is better than inequality \eqref{maxiteratedSob}.}
This follows from the fact that $M^kf\leq c_{r,k}\, M_{r}f$ since $M_rf \in A_1$ for $r>1$, with constant independent of $f$.

Inequality \eqref{Sobimpro} holds in the range $1<p<n$, but it is false when $p=1$. 
{In this case we have two different replacements. First, we  claim that}
\begin{equation*} \label{GNSimprweakI}
\|M_{n'}f\|_{L^{n',\infty}(w^{n'})}\leq c\, {[w]_{A_{1,n'}}^{\frac{2}{n'}} } \|w\nabla f\|_{{L^1(\R^n)}},
\end{equation*}
using the classical $A_1$ weak type $(1,1)$ estimate of Fefferman-Stein,
$$
\left \Vert Mf\right \Vert_{L^{1,\infty}(u)} \leq c\,\Vert f \Vert_{L^{1}(Mu)},
$$
combined with \eqref{oneweightSob} when $p=1$. 

Alternatively, consider the quantity $\|wM_{n'}f\|_{L^{n',\infty}(\R^n)}$ which is very different from $\|M_{n'}f\|_{L^{n',\infty}(w^{n'})}$.  We claim that for a constant $c$ depending on the dimension and $ [w]_{A_{1,n'}}$ we have
\begin{equation*} \label{GNSimprweakII}
\|wM_{n'}f\|_{L^{n',\infty}(\R^n)}\leq c\, \|w\nabla f\|_{{L^1(\R^n)}}.
\end{equation*}
The difficulty is that the weight is inside the distribution level of the definition of the norm making the problem non-standard. 
However, it is proved in \cite{CMP} (improved in \cite{LiOP} and then in \cite{PR})  that the Hardy-Littlewood maximal operator satisfies a Sawyer-type estimate of the form
\begin{equation} \label{cmpbd}
\left \Vert \frac{Mf}{v}\right \Vert_{L^{1,\infty}(uv)} \leq c_{[u]_{A_1},[uv]_{A_{\infty}}} \Vert f \Vert_{L^{1}(u)},
\end{equation}
where $u\in A_1$ and $uv\in A_{\infty}$.

Applying this  with $v=w^{-n'}$ and $u=\frac1v$, since $w^{n'} \in A_1$ we have
%
%
%
%
\begin{align*}
\|wM_{n'}f\|_{L^{n',\infty}(\R^n)} & =\|w^{n'}M(|f|^{n'})\|_{L^{1,\infty}(\R^n)}^{\frac1{n'}} \\
& \leq c_{ [w]_{A_{1,n'}}} \, \|w^{n'}|f|^{n'}\|_{L^{1}(\R^n)}^{\frac1{n'}} \\
& =  c_{ [w]_{A_{1,n'}}} \,  \|wf\|_{L^{n'}} \leq c_{[w]_{A_{1,n'}}}\,\|w\nabla f\|_{{L^1(\R^n)}}
\end{align*}
for which we have to use the inequalities \eqref{oneweightSob} and \eqref{cmpbd}.  The explicit dependence on the constant is not ideal and we do not state it here.

Interestingly enough, inequality \eqref{PointMn'} is false for the pointwise larger operator $M_{n'+\epsilon}$, when $\epsilon>0$. It is possible to replace $M_{n'}$ by a slightly larger maximal operator defined in the scale of the  Lorentz space $L^{n',1}$.   Namely, it is shown in \cite{HMP}, that 
\begin{equation}\label{PointMn',1}
M_{n',1}f(x) \leq c\,I_{1}(|\nabla f|)(x),
\end{equation}
which yields an improvement of \eqref{Sobimpro} and hence of \eqref{oneweightSob} when $p\in(1,n)$. Here we are using the normalized definition of the 
\[
M_{ s,q }f(x) =
\sup_{x \in Q}
\frac{1}{ |Q |^{1/s} } \norm{  f\, \chi_{_{Q}} }_{ {_{L^{s,q}(dx)} } },
\]
$0 < s,q \leq \infty$, where
\[ \norm{f}_{_{L^{s,q}(\mu)}} = \left[ s\, \int_{0}^{\infty}  \left(
t\, \mu\{  x \in \R^{n}: \absval{f(x)} > t \}^{1/s} \right)^{q}\,
\frac{dt}{t} \right]^{1/q} < \infty,
\]
{if  $q < \infty$, and \,$\sup_{ 0
< t < \infty}{t\, \mu\{  x \in \R^{n}: \absval{f(x)} > t \}^{1/s}} < \infty,
$}
when $q = \infty$.  The proof of inequality \eqref{PointMn',1} is based on the following improvement of the classical Poincar\'e-Sobolev inequality in Lorentz spaces due to O'Neil \cite{ON} and Peetre \cite{Pee} (see also \cite{MP}), 
\begin{equation*}\label{rspoin}\|f-f_Q\|_{L^{p^*,p}(Q)}\leq c\, \ell(Q)\left(\dashint_Q|\nabla f|^p\right)^{\frac1p},\end{equation*}
for $1\leq p<n$, where we are using the normalized Lorentz average.  In particular, when $p=1$, we have the following 
\begin{equation*}\label{lorentzPoin}\|f-f_Q\|_{L^{n',1}(Q)}\leq c\,  \ell(Q)\, \dashint_Q |\nabla f|\end{equation*}
from which we can derive  \eqref{PointMn',1} easily. 

A natural question would be to find the largest maximal function $M_\Phi$  for which \eqref{PointMn',1} holds, where 
$\Phi$ is a Young function and $M_\Phi$ is defined through a Luxemburg type norm.

%

\subsection{The case of singular integrals} \label{SIO}

Motivated by maximal operator results discussed above, we proved similar results for more complicated singular integral operators in \cite{HMP}. This will be  the case of several classical smooth convolution singular integrals but we will be considering operators beyond Calder\'on-Zygmund operators like the so-called rough singular integral operators $T_{\Omega}$.  Surprisingly, we found Sobolev type estimate as \eqref{maxSob} or \eqref{Sobimpro} for these operators which are defined by the expression  
$$T_\Omega f(x)=p.v.\int_{\R^n}\frac{\Omega(y')}{|y|^n}f(x-y)\,dy$$
where $y'=\frac{y}{|y|} \in \mathbb S^{n-1}$ which are well defined if $\Omega\in L^1(\mathbb S^{n-1})$ and has zero average on $\mathbb S^{n-1}$.  We are not able to work with the assumption $\Omega\in L^1(\mathbb S^{n-1})$ to obtain our results, but rather $\Omega \in L^{n,\infty}(\mathbb S^{n-1})$ which properly contains the spaces $L^r(\mathbb S^{n-1})$ for $r\geq n$.

%
%
%
%
%
%

One of our main results from  \cite{HMP} is the following.

\begin{theorem}[\hspace{-.2mm}\cite{HMP}] \label{rough}  Let $T_{\Omega}$ be a rough singular integral convolution operator as above with $\Omega \in L^{n,\infty}(\mathbb S^{n-1})$.  Then
\begin{equation}\label{ptwisebdd} 
|T_{\Omega} f(x)|\leq c\, \|\Omega\|_{L^{n,\infty}(\mathbb S^{n-1})}\, I_1(|\nabla f|)(x), \quad f\in C^\infty_c(\R^n).
\end{equation}
Hence, if $p\in (1,n)$ and if we let $w\in A_{p,p^*}$, we have
\begin{equation}\label{SIOroughStrong}
\|wT_{\Omega}f\|_{L^{p^*}(\R^n)}\leq c_{p}\, \|\Omega\|_{L^{n,\infty}(\mathbb S^{n-1})}\,
 [w]_{A_{p,p^*}}^{\frac1{n'}\max\{1,\frac{p'}{p^*}\}}\, \|w\nabla f\|_{L^p(\R^n)} 
\end{equation}
and 
\begin{equation*}\label{SIOroughWeak}
\|T_{\Omega}f\|_{L^{n',\infty}(w^{n'})}\leq  c\, \|\Omega\|_{L^{n,\infty}(\mathbb S^{n-1})}\,[w]_{A_{1,n'}}^{\frac1{n'}}\,   \|w\nabla f\|_{L^1(\R^n)}.
\end{equation*}

\end{theorem}
{Similar results hold for rough maximal} type operators and more exotic maximal operators such as the spherical maximal operator.

We note here that all the singular integral operators for which we have proven \eqref{ptwisebdd} have a cancellation on the kernel. 
We do not know} if \eqref{ptwisebdd} also holds for Calder\'on-Zygmund operators of non-convolution type.
We do however have reasonable bounds for these operators.
Recall that these are a continuous linear operator $T:C_0^{\infty}(\R^n) \to
\mathcal{D}'(\R^n)$ { which  extend
to bounded operators on} $L^2(\R^n)$, and whose distributional
kernel coincides away from the diagonal $x=y$ in $\R^n\times
\R^n$, with a function $K$ so that
$$Tf(x)=\int_{{\mathbb R}^n}K(x,y)f(y)dy, \quad f\in C_0^{\infty}(\R^n), x\not\in\mbox{supp}(f),
$$
and satisfies the standard estimates, namely, the size estimate
\begin{equation}\label{cond1}
|K(x,y)|\le \frac{c}{|x-y|^n}
\end{equation}
and the regularity condition: for some $\epsilon>0$,
\begin{equation}\label{cond2}
|K(x,y)-K(z,y)|+|K(y,x)-K(y,z)|\le c\, \frac{|x-z|^{\epsilon}}{|x-y|^{n+\epsilon}},
\end{equation}
whenever $2\,|x-z|<|x-y|$. 

By the {well-known} classical theory, we have 
\begin{equation}\label{scztype}
T: L^{p}(\R^n) \to L^{p}(\R^n)
\end{equation}
 when $1<p<\infty$ and
\begin{equation}\label{wcztype}
T: L^{1}(\R^n) \to L^{1,\infty}(\R^n).
\end{equation}

We will be using the following result due to  Coifman-Fefferman \cite{CF}, 
\begin{equation}\label{CF}
\|T f\|_{L^p(w)}\leq c_{p,T}\,\,[w]_{A_{\infty}} \, \|Mf\|_{L^p(w)}, 
\end{equation}
for $p\in (0,\infty)$ and $w\in A_{\infty}$. The original proof does not yield this precise bound since some further ideas are required. See details in \cite{OPR}. A similar result holds for the weak norm.

Now, arguing as in Section \ref{ToyExample}, {using the sharp weighted estimate for Calder\'on-Zygmund operators,
%
%
since $w^{p^*} \in  A_{p^*} $, we have}
%
%
%
\begin{align*}
\|wTf\|_{L^{p^*}(\R^n)} & \leq c_{p,T}\,[w^{p^*}]_{A_{p^*}}^{\max\{  1, \frac{1}{p^*-1}  \}} \, \|wf\|_{L^{p^*}(\R^n)} \\
& \leq c_{p,T}\,[w]_{A_{p,p^*}}^{\max\{  1, \frac{1}{p^*-1}  \}+\frac1{n'}} \,[w]_{A_{p,p^*}}^{\frac1{n'}}\,\|w\nabla f\|_{L^{p}(\R^n)}
\end{align*}
in the case $p\in(1,n)$, which is certainly less optimal than \eqref{SIOroughStrong}.  Moreover, we may use the fact that $A_{\frac{p^*}{n'}}$ is a proper subclass of $A_{p^*}$ and the estimate of Duoandikoetxea \cite{Duo} (see also \cite{LM})
\begin{equation*} 
\|Tf\|_{L^p(w)}\leq c_{p,q,T}\,[w]_{A_q}\|f\|_{L^p(w)},\quad 1\leq q<p,
\end{equation*}
to obtain when $p\in (1,n)$
\begin{equation} \label{DuoCZO}
\|wTf\|_{L^{p^*}(\R^n)}\leq c_{p,T}\,[w]_{A_{p,p^*}}^{1+\frac1{n'}}\|w\nabla f\|_{L^p(\R^n)}.
\end{equation}
On the other hand, if we use \eqref{CF} instead, we claim 
%
\begin{equation} \label{claimTCZO}
\|wTf\|_{L^{p^*}(\R^n)} \leq 
c_{p,T}\, [w^{p^*}]_{A_{\infty}}  \, [w]_{A_{p,p^*}}^{\frac1{n'}\max\{1,\frac{p'}{p^*}\}} \, \|w\nabla f\|_{L^p(\R^n)}.
\end{equation}
Indeed, 
\begin{align*}
\|wTf\|_{L^{p^*}(\R^n)} &= \|T f\|_{L^{p^*}(w^{p^*}) } 
  \leq c_{p,T}\,  [w^{p^*}]_{A_{\infty}}\, \|M f\|_{L^{p^*}(w^{p^*}) } \\
&\leq 
c_{p,T}\, [w^{p^*}]_{A_{\infty}}\, [w]_{A_{p,p^*}}^{\frac1{n'}\max\{1,\frac{p'}{p^*}\}} \, \|w\nabla f\|_{L^p(\R^n)}
\end{align*}
by \eqref{maxSob},  which is also a worse  result than  \eqref{SIOroughStrong}.

Finally, using the main result in \cite{LOP},  we also have 
\begin{align*}
\|T f\|_{L^{n',\infty}(w^{n'})} &\leq 
c_{T}\,[w^{n'}]_{A_1}\, \|f\|_{L^{n'}(w^{n'})}\\
& \leq c_{T}
\, [w]_{A_{1,n'}}^{1+\frac1{n'}}\,   \|w\nabla f\|_{L^1(\R^n)},
\end{align*}
which is a better result than \eqref{DuoCZO}.

\section{The linear commutators} \label{linear commutators}

Likely motivated by the work of A. P. Calder\'on on commutators, R. Coifman, R. Rochberg and G. Weiss introduced in \cite{CRW} the operator
\begin{equation*}
T_b f(x)=\int_{\R^n} (b(x)-b(y))K(x,y)f(y)\,dy,
\end{equation*}
where $K$ is a kernel satisfying the standard Calder\'on-Zygmund
estimates (see Section \ref{SIO}) and where $b$ is any locally integrable function. The
operator is called commutator since 
$$T_b= [b,T]=b T -T(b\,\cdot)$$ 
where $T$ is the Calder\'on-Zygmund singular
integral operator associated to $K$. The main result
from \cite{CRW} states that $[b, T]$ is a bounded operator on
$L^p(\R^n)$ { if} $1<p<\infty$ and $b$ is a
$BMO$ function. They proved that
\begin{equation*}\label{bmolptolp}
\|T_b\|_{L^p(\R^n) \to L^p(\R^n)} \leq  c_{p,T}\, \|b\|_{BMO}.
\end{equation*}
Moreover, they provided a new and beautiful characterization of $BMO$ in the case of $T$ being the Hilbert transform or any Riesz transform $R_j$, i.e.
\begin{equation*}
\|b\|_{BMO} \approx   \|T_b\|_{L^p(\R^n) \to L^p(\R^n)}.
\end{equation*}

These commutators have proved to be of great interest in many
situations, for instance, the theory of non divergence elliptic equations with discontinuous
coefficients \cite{CFL1} \cite{CFL2} \cite{DiR}. We refer to \cite{I} for a beautiful account.

A natural generalization of the commutator $[b,T]$ is given by
\begin{equation*}\label{conmutador}
T^m_b f(x)=\int_{\mathbb R^n} (b(x)-b(y))^m K(x,y)f(y)\,dy,
\end{equation*}
where $m\in \mathbb N$, and the case $m=0$ recaptures the
Calder\'on-Zygmund singular integral operator. Observing that for $m\geq 2$, the iterated commutator
$$T^m_b = \overbrace{[b,\cdots, [b,T]] }^{(m\, times)}$$
is also a bounded operator 
$$T^m_b: L^p(\R^n) \to L^p(\R^n),$$
whenever $1<p<\infty$ and $b$ is a $BMO$ function.

It was shown in \cite{P1} that these operators are not of weak type $(1,1)$ and that there is an intimate relationship with the iterated maximal function $M^m$ We will need the following result: for any $0<p<\infty$ and
any $w\in A_{\infty}$, we have
%
\begin{equation}\label{c1}
\|T^m_b f\|_{L^p(w)}\leq c_{m,p,T}\,[w]_{ A_{\infty}}   \,\|b\|_{BMO}^{m} \|M^{m+1}(f)\|_{L^p(w)}.
\end{equation}
This inequality is sharp since  $M^{m+1}$ can not be
replaced by the smaller operator $M^m$.   $[w]_{A_{\infty}}$ denotes the usual  $A_{\infty}$ constant, 
$$
[w]_{A_\infty}:=\sup_Q\frac{1}{w(Q)}\int_Q M(w\chi_{Q})\,dx.
$$
%

On the other hand, it follows from \cite{CPP} that 
\begin{align*}
\|wT^m_b f\|_{L^{p^*}(\R^n)} & = \|T^m_b f\|_{L^{p^*}(w^{p^*})} 
\\
& \leq c_{m,p,T}\,[w^{p^*}]_{A_{p^*}}^{(m+1)\max\{1,\frac1{p^*-1}\}} \|b\|_{BMO}^{m} \|f\|_{L^{p^*}(w^{p^*})},
\end{align*}
which combined with   \eqref{oneweightSob} yields
$$
\|wT^m_b f\|_{L^{p^*}(\R^n)} \leq c_{m,p,T}\, [w]_{A_{p,p^*}}^{(m+1)\max\{1,\frac1{p^*-1}\}} \|b\|_{BMO}^{m} [w]_{A_{p,{p^*}}}^{\frac1{n'}}\|w\nabla f\|_{L^p(\R^n)}.
$$
The above constant is worse than the one derived using  \eqref{c1}  $[w^{p^*}]_{A_{\infty}}$ times $[w]_{A_{p,{p^*}}}^{\frac1{n'}}$) at least in the case $p$ is close to $1$ and $n$ is large.

Reasoning as in the proof of \eqref{claimTCZO}, for { $p\in (1,n)$} we have 
\begin{equation*} \label{claimCommut}
\|wT^m_bf\|_{L^{p^*}(\R^n)} \leq 
c_{m,p,T}\, [w^{p^*}]_{A_{\infty}}^{m+1}\,[w]_{A_{p,p^*}}^{\frac1{n'}\max\{1,\frac{p'}{p^*}\}} \, \|w\nabla f\|_{L^p(\R^n)}.
\end{equation*}
The case $p=1$  has a different nature, observe that the exponent above blows up when  $p\to1$.  We end this section with the following result whose proof is based on the $A_1$ results for commutators by C. Ortiz in \cite{O}.
%
\begin{theorem} \label{CommSob} {  Let $w\in A_{1,n'}$. Then}

\begin{equation*}
\|T^m_b f\|_{L^{n',\infty}(w^{n'})} \leq c_T\,\|b\|_{BMO}^{m}\,[w]_{A_{1,n'}}^{m+1+\frac1{n'}} \|w\nabla f\|_{L^1(\R^n)}.
\end{equation*}
\end{theorem}

\section{The nonlinear commutators}\label{nonlinear commutators}

The second kind of commutators that we are going to consider were introduced by R.
Rochberg and G. Weiss in \cite{RW}. These are nonlinear operators defined for
appropriate functions by
\[
f
{\mapsto}
Nf= N_T(f):=T(f\,\log \absval{f} )  - Tf\, \log \absval{Tf}.
\]
%
%
%
We will assume that $T$ is a Calder\'on-Zygmund singular integral operator in this section.
It is of interest to study this operator due to its
relationship with the Jacobian mapping and with nonlinear P.D.E.
as shown in \cite{IS} \cite{GI} (see also \cite{M}).

The main result from \cite{RW} is the following $L^p$ bound
\begin{equation*}
\norm{
Nf}_{L^{p} (\R^{n})} \le c_{p,T}\, \norm{f }_{L^{p} (\R^{n})}, \quad 1<p<\infty.
\end{equation*}
We have the following improvement from \cite{P2} that uses a completely different method. 

\begin{theorem} \label{CFNLC}
Suppose that $0 < p < \infty$ and $w\in A_{\infty}$. Then,
\begin{equation*}
\|Nf\|_{L^p(w)}\leq  c_{p,T}[w]_{A_{\infty}}\, \log\left( e+[w]_{A_{\infty}} \right)\, \|M^{2}f\|_{L^p(w)}.
\end{equation*}
\end{theorem}
We briefly sketch the proof below.  The idea is to split the nonlinear commutator $N$ {into simpler pieces:}
\begin{align*}
T(f\,\log \absval{f} )  - Tf\, \log \absval{Tf}
&=T(f\,\log \frac{\absval{f}}{Mf} ) +
[\log Mf,T]f  - Tf\, \log \frac{ \absval{Tf}}{Mf} \\
&:=N_{1}f + N_{2}f + N_{3}f.
\end{align*}

%
To estimate $N_{1}$, we write
\[
N_{1}f=T\left(Mf\,
\frac{f}{Mf}\,  \log ( \frac{ \absval{f} }{Mf} )
\right),
\]
and we apply \eqref{c1} with $m=0$. Since $w\in A_{\infty}$, we have
\begin{align*}
\int_{\R^{n}} \absval{ N_{1} f}^{p}\, w
& \le c\, [w]_{A_{\infty}}^{p}\, \int_{\R^{n}}
\left(M\left(Mf\, \frac{f}{Mf}\,  \log (  \frac{ \absval{f} }{Mf} ) \right)
\right)^{p} \, w \\
& \le c\,
[w]_{A_{\infty}}^{p}\, \int_{\R^{n}}
(M^{2}f)^{p}\,
w
\end{align*}
by utilizing the fact that $\absval{t\,\log t} \le  \frac{1}{e}$ when $0<t\le 1$.

For $N_{2}$, we notice that it is of commutator type fixing $b_f= \log Mf \in BMO$  with dimensional constant, namely  independent of $f$ and hence we can apply  again  \eqref{c1} (see \cite{P2} for details)
\begin{align*}
\int_{\R^{n}} \absval{ N_{2}f }^{p}\, w
&=
\int_{\R^{n}} \absval{ [\log Mf,T]f }^{p}\, w \\
&\le c\, \norm{\log Mf }_{BMO }^{2p} \, [w]_{A_{\infty}}^{p}
\, \int_{\R^{n}} (M^{ 2 }f)^{p}\,  w\\
&\le
c\, [w]_{A_{\infty}}^{p}\, \int_{\R^{n}} (M^{ 2 }f)^{p}\,  w.
\end{align*}

For $N_{3}$, we split $\R^{n}$ in two disjoint sets $A$ and
$B$ where 
$$A= \{y \in \R^{n} :\absval{Tf(y)} \le Mf(y)  \}$$ and $B=  \R^{n}\setminus A.$  Writing
\[
N_{3}f = Mf\, \frac{ \absval{Tf} }{Mf} \, \log \frac{ \absval{Tf} }{Mf}
\]
we have that $|N_3f|\leq Mf$ on $A$, and on $B$ we use $\log t \leq\, \frac{  t^{\epsilon} 
}{\epsilon}$, $t>1$, $\epsilon >0$ to get
$$
\|N_3f\|_{L^p(w)}\leq  \|N_3f \chi_A\|_{L^p(w)}+\|N_3f \chi_B\|_{L^p(w)} 
$$
$$
\leq  \|Mf \|_{L^p(w)}+ \frac1{\epsilon}\,
\| \absval{Tf}^{(\epsilon +1)}\,(Mf)^{ - \epsilon} \|_{L^p(w)} 
$$
$$
=  \|Mf \|_{L^p(w)}+ \frac1{\epsilon}\,
\| Tf\|_{L^{p(1+\epsilon)}(w(Mf)^{ - \epsilon p})}^{1+\epsilon} 
$$

%
To apply \eqref{c1} again with $m=0$, we must verify that \, $(Mf)^{ - \epsilon \, p}\,w \in
A_{\infty}$\, for $\epsilon$ small enough and with a constant
independent of $f$. Being more precise, it can be shown that
$$[wu^{-1}]_{ A_{\infty}}\leq [w]_{A_{\infty}} \,[u]_{A_1},
$$
if $\epsilon$ is chosen to be $\epsilon p<1 $. Indeed,  we have by the well-known Coifman-Rochberg theorem that if\, $0<\delta<1$,\, 
$u=(Mf)^{1-\delta}\in A_1$ with  $[(Mf)^{1-\delta}]_{A_1} \leq \frac{c_n}{\delta}$ and then
$$[w(Mf)^{ - \epsilon \, p}]_{ A_{\infty}}\leq \frac{c_{n}}{1-\epsilon p}\, [w]_{A_{\infty}},
$$
from which  we have
$$
\|N_3f\|_{L^p(w)}\leq \|Mf \|_{L^p(w)}+ 
\frac1{\epsilon}\, [w(Mf)^{ - \epsilon p}]_{A_{\infty}}^{1+\epsilon}c_{p,T}^{1+\epsilon}
\| Mf\|_{L^{p(1+\epsilon)}(w(Mf)^{ - \epsilon p})}^{1+\epsilon}
$$
$$
\leq \|Mf \|_{L^p(w)}+  \frac1{\epsilon}\, [w]_{A_{\infty}}^{1+\epsilon}\, [(Mf)^{  \epsilon p}]_{A_{1}}^{1+\epsilon}\, c_{p,T}^{1+\epsilon}
\| Mf\|_{L^{p}(w)}
$$
$$
\leq \|Mf \|_{L^p(w)}+  \frac1{\epsilon}\, [w]_{A_{\infty}}^{1+\epsilon}\, (\frac{c_n}{1-p\epsilon})^{1+\epsilon}
\,c_{p,T}^{1+\epsilon}
\| Mf\|_{L^{p}(w)}
$$
Optimizing this estimate  with $\epsilon=\frac{1}{p \log (e^2 + [w]_{A_{\infty}  }) }$ we get finally
$$
\|N_3f\|_{L^p(w)}\leq  [w]_{A_{\infty}} \log\left( e+[w]_{A_{\infty}} \right)\,c_{p,T}\,
\|Mf \|_{L^p(w)}
$$
and hence 
\begin{equation} \label{corNonlinear}
\|Nf\|_{L^p(w)}\leq c_{p,T}\,[w]_{A_{\infty}} \log\left( e+[w]_{A_{\infty}} \right)\, \|M^2f \|_{L^p(w)}
\end{equation}
%

This yields the proof of Theorem \ref{CFNLC}.

Using \eqref{BuckBd} we have as an immediate consequence the following corollary.

\begin{corollary}
Let $1 < p < \infty$ and $w\in A_{p}$. Then,
\begin{equation*}
\|Nf\|_{L^p(w)}\leq c_{p,T}\,[w]_{A_{\infty}} \log\left( e+[w]_{A_{\infty}} \right)\, [w]_{A_{p}}^\frac{2}{p-1} \|f\|_{L^p(w)}.
\end{equation*}
\end{corollary}

Another result that can { be} found in \cite{P2} is the following 

\begin{theorem}  \label{principalnl}
Suppose that $1 < p < \infty$. {Let $w$ be any weight.} Then
\begin{equation*}
\int_{\R^{n}} \absval{ Nf(x) }^{p}\, w(x)\,dx
\le c\int_{\R^{n}} \absval{f(x)}^{p}\,  M^{ [2 p ]+1  }w(x)\,dx.
\label{NLC}
\end{equation*}
\end{theorem}

The method used in \cite{P1,P2} cannot be applied here since the operator is non-linear. The splitting above has to be used.

Finally, combining   \eqref{corNonlinear} and the Sobolev type inequality \eqref{maxSobIter}  derived in \cite{HMP}
\begin{equation*} \label{claimNLC}
\|wNf\|_{L^{p^*}(\R^n)} \leq 
c_{p,T}\, [w^{p^*}]_{A_{\infty}}\,\log\left( e+[w^{p^*}]_{A_{\infty}} \right) \,[w]_{A_{p,p^*}}^{\frac1{n'}\max\{1,\frac{p'}{p^*}\}}\,\|w\nabla f\|_{L^p(\R^n)}.
\end{equation*}
Observe that, as in the linear case, the exponent blows up when  $p\to1$.  
We state now a conjecture motivated by Theorem \ref{CommSob} without such a blow up.

\begin{conjecture} \label{ConjecN} {  Let $w\in A_{1,n'}$. Then}
\begin{equation*}
\|Nf\|_{L^{n',\infty}(w^{n'})} \leq c\, 
[w]_{A_{1,n'}}^{2+\frac1{n'}}\,
\|w\nabla f\|_{L^1(\R^n)}.
\end{equation*}
\end{conjecture}


\begin{thebibliography}{99}




\bibitem[Buc]{Buc} {\sc S. Buckley}, \emph{Estimates for operator norms on weighted spaces and reverse Jensen inequalities}. Trans. Amer. Math. Soc. {\bf 340} (1993), 253--272.







\bibitem[CF1]{CFL1}
{\sc F. Chiarenza}, {\sc M. Frasca,} and {\sc P. Longo}, {\em Interior $W^{2,p}$
estimates for non divergence elliptic equations with discontinuous
coefficients}, Richerche Mat. {\bf 40} (1991), 149--168.

\bibitem[CF2]{CFL2}
{\sc F. Chiarenza}, {\sc M. Frasca,} and {\sc P. Longo}, {\em $W^{2,p}$--solvability of
the Dirichlet problem for nondivergence elliptic equations with
VMO coefficients}, Trans. Amer. Math. Soc. {\bf 334} (1993),
841--853.


\bibitem[CPP]{CPP} {\sc D. Chung}, {\sc C. Pereyra} and {\sc C. P\'erez}, Sharp bounds for
general commutators on weighted Lebesgue spaces, Trans. Amer. Math.
Soc., 364 (2012), 1163\textendash 1177.


\bibitem[CF]{CF} {\sc R. Coifman} and {\sc C. Fefferman}, \emph{Weighted norm inequalities for maximal functions and singular integrals}, Studia Math. 51 (1974), 241--250.





\bibitem[CRW]{CRW} {\sc R. R. Coifman}, {\sc R. Rochberg} and {\sc G. Weiss}, {\it Factorization theorems for Hardy spaces in several variables}, Ann. of Math. (2) 103 (1976), no. 3, 611\textendash635.

\bibitem[CMP]{CMP} {\sc D. Cruz-Uribe}, {\sc J.M. Martell}, and {\sc C. P\'erez}, {\em
Weighted weak-type inequalities and a conjecture of Sawyer,}
{Int. Math. Res. Not. IMRN,} {\bf 30} 2005, 1849-1871.

\bibitem[Dir]{DiR} {\sc G. Di Fazio} and {\sc M. A. Ragusa}, \emph{Interior estimates in Morrey spaces
for strong solutions to nondivergence form equations with
discontinuous coefficients}, { J. Funct. Anal.}  {\bf 112} (1993), 241--256.


\bibitem[Duo]{Duo} {\sc J. Duoandikoetxea}, \emph{Extrapolation of weights revisited: new proofs and sharp bounds}. J. Funct. Anal. 260 (2011), 1886--1901.



\bibitem[GI]{GI}
{\sc L. Greco,} and {\sc T. Iwaniec} {\em New inequalities for the Jacobian},
Ann. Inst. Henri Poincar\'e, {\bf 11} (1994), 17--35.

\bibitem[Haj]{Haj} {\sc P. Haj\st{l}asz}, \emph{Sobolev inequalities, truncation method, and John domains}. Papers on Analysis, Rep. Univ. Jyv\"askyl\"a Dep. Math. Stat. 83, 2001, 109--126.

\bibitem[HMS]{HMS} {\sc E. Harboure}, {\sc R. Mac\'ias}, and {\sc C. Segovia}, \emph{Extrapolation results for classes of weights}, Amer. J. Math., {\bf 110} (1988), 383-397.


\bibitem[HMP]{HMP}  {\sc C. Hoang}, {\sc K. Moen}, and {\sc C. P\'erez}, 
{\em Pointwise estimates for rough operators with applications to weighted Sobolev inequalities,} to appear Journal d'Analyse Mathématique .


\bibitem[I]{I}  {\sc T. Iwaniec}, {\it Let the Beauty of Harmonic Analysis be Revealed Through Nonlinear PDEs; A Work of Art in Three Sketches},  Nov 2012; Rev. Mat. Iberoam., Special Issue. 

\bibitem[IS]{IS}
{\sc T. Iwaniec} and {\sc C. Sbordone}, {\em Weak minima of variational integrals},
J. Reine Angew Math. {\bf 454}, 143--161.

\bibitem[KLV]{KLV} {\sc J. Kinnunen, J. Lehrb\"ack}, and {\sc A. V\"ah\"akangas}, \emph{Maximal function methods for Sobolev spaces}, Mathematical Surveys and Monographs, Volume: 257, 2021.

\bibitem[LMPT]{LMPT} {\sc M.T. Lacey, K. Moen, C. P\'erez}, and {\sc R.H. Torres}, \emph{Sharp weighted bounds for fractional integral operators.} J. Funct. Anal. 259 (2010), no. 5, 1073--1097.

\bibitem[LM]{LM} {\sc A. Lerner} and {\sc K. Moen}, \emph{Mixed $A_p$-$A_\infty$ estimates with one supremum.}
Studia Math. 219 (2013), no. 3, 247--267. 

\bibitem[LOP]{LOP} {\sc A. Lerner, S. Ombrosi} and {\sc C. P\'erez}, {\em Sharp $A_1$ bounds for
Calder\'on-Zygmund operators and the relationship with a problem of
Muckenhoupt and Wheeden,} 
{Int. Math. Res. Not. IMRN,} 2008,  \textbf{6}, Art. ID rnm161, 11 pp. 42B20.


\bibitem[LiOP]{LiOP} {\sc K. Li, S. Ombrosi}, and {\sc C. P\'erez}, \emph{Proof of an extension of E. Sawyer's conjecture about weighted mixed weak-type estimates.} {Math. Ann.} (2018).

\bibitem[LPR]{LPR} {\sc T. Luque, C. P\'erez} and {\sc E. Rela}, {\em Reverse H\"older property for strong weights and general measures,}  J. Geom. Anal. 27 (2017), no. 1, 162--182. 

\bibitem[MP]{MP} {\sc J. Mal\'y} and {\sc L. Pick}, \emph{An elementary proof of sharp Sobolev embeddings}. Proc. Amer. Math. Soc. {\bf 130} (2002), no. 2, 555--563.


\bibitem[Maz]{Maz} {\sc V. Maz'ya}, \emph{Sobolev spaces with applications to elliptic partial differential equations.} Second, revised and augmented edition, Springer-Verlag Berlin Heidelberg (2011).

\bibitem[M]{M} {\sc M. Milman}, Extrapolation and Optimal Decompositions, Lect. Notes Math. 1580, Springer Verlag, (1995).

\bibitem[MW]{MW} {\sc B. Muckenhoupt} and {\sc R. Wheeden}, {\em Weighted norm inequalities for fractional integrals}, Trans. Amer. Math. Soc. 192 (1974) 261-274.


\bibitem[ON]{ON} {\sc R. O'Neil}, {\em Convolution operators and L(p,q) spaces}. Duke Math. J. {\bf 30} 1963, 129--142.


\bibitem[O]{O} {\sc C. Ortiz-Caraballo}, {\em Quadratic $A_{1}$ bounds for commutators
of singular integrals with $BMO$ functions,} Indiana Univ. Math. J. 60
(2011), no. 6, 2107\textendash 2130.

\bibitem[OPR]{OPR} {\sc C. Ortiz-Caraballo},  {\sc C. P\'erez} and {\sc E. Rela},  {\em Improving bounds for singular operators via Sharp Reverse H\"older Inequality for $A_{\infty}$, }   ``Operator Theory: Advances and Applications",  {\bf 229}, (2013),  303-321. Springer Basel. 




\bibitem[Pee]{Pee} {\sc J. Peetre}, 
{\em 
Espaces d'interpolation et    { th\'eor\`em} de Soboleff,
} 
Ann. Inst. Fourier {\bf 16} 1966, 279-317.




\bibitem[P1]{P1}  {\sc C. P\'erez}, {\em Endpoint estimates for
commutators of singular integral operators}, { J. Funct. Anal.} (1) {\bf 128} (1995), 163--185.




\bibitem[P2]{P2} {\sc C. P\'erez}, {\em Sharp estimates for commutators of singular
integrals via iterations of the Hardy--Littlewood maximal function},
{J. Fourier Anal. Appl.} (6) {\bf 3} (1997),
108--146.

\bibitem[PR]{PR}  {\sc C. P\'erez} and {\sc E. Roure}, {\em Sawyer-type inequalities for Lorentz spaces,}  
Math. Ann. 383 (2022), no. 1-2, 493--528.

\bibitem[RW]{RW} 
{\sc R. Rochberg} and {\sc G. Weiss}, {\it  Derivatives of analytic families of Banach spaces,} Ann. of Math. (2) 118 (1983), no. 2, 315--347. 

\bibitem[Sal]{Sal}  {\sc L. Saloff-Coste}, \emph{Aspects of Sobolev-Type Inequalities}, (London Mathematical Society Lecture Note Series) (2001).





\bibitem[St]{St} {\sc E. Stein}, {\em Singular integrals and differentiability properties of functions,} Princeton Mathematical Series, No. 30 Princeton University Press, Princeton, N.J.

\bibitem[StW]{StW} {\sc E. Stein} and {\sc G. Weiss}, {\em Fractional integrals on n-dimensional Euclidean space,} J. Math. Mech., 7 (1958),
503-514.


\end{thebibliography}
\end{document}